\documentclass{amsart}

\newcommand\neumanncite[1]{\cite{#1}}

\usepackage{xypic}
\newbox\neumannbox
\def\overNeTag#1#2#3{\setbox\neumannbox\hbox{$#1$}\hbox to
  0pt{\vbox to 0pt{\vglue-#3\vglue-\ht\neumannbox\hbox to \wd\neumannbox
      {\hss$\scriptstyle#2$\hss}\vss}\hss}\box\neumannbox}
\def\underNeTag#1#2#3{\setbox\neumannbox\hbox{$#1$}\hbox to 0pt{\vbox to
    0pt{\vglue#3\vglue\ht\neumannbox\hbox to \wd\neumannbox
      {\hss$\scriptstyle#2$\hss}\vss}\hss}\box\neumannbox}
\def\leftNeTag#1#2#3{\hbox to 0pt{\vbox to 0pt{\vss\hbox to
      0pt{\hss$\scriptstyle#2$\hskip#3}\vss}}#1}
\def\rightNeTag#1#2#3{\hbox to 0pt{\vbox to 0pt{\vss\hbox to
      0pt{\hskip#3$\scriptstyle#2$\hss}\vss}}#1}
\def\CirC{\lower.2pc\hbox to 2pt{\hss$\circ$\hss}}

\newtheorem{theorem}{Theorem}[section]
\newtheorem{proposition}{Proposition}[section]
\newtheorem{conjecture}{Conjecture}[section]
\theoremstyle{definition}

\newtheorem{question}{Question}[section]

\begin{document}

%\title{Principal analytic link-theory in surface singularity links} 
% Capitaliztion should be done by the style file!!!!
%\title{PRINCIPAL ANALYTIC LINK-THEORY IN SURFACE SINGULARITY LINKS} 
\title{Complex analytic realization of links}
\author{Walter D. NEUMANN}

\address{Department of Mathematics, Barnard
  College, Columbia University\\ New-York, NY 10027, USA \\
E-mail: neumann@math.columbia.edu}

\author{Anne  PICHON}
\address{Institut de Math\'ematiques de Luminy\\ UMR 6206 CNRS,
Campus de Luminy - Case 907\\ 13288 Marseille Cedex 9,
France\\
E-mail: pichon@iml.univ-mrs.fr}

%\noindent
%{\bf AMS Subject Classification :} 32S50, 14B05, 57M25
%\vspace{.1in}
%\noindent
\begin{abstract}
  We present the complex analytic and principal complex analytic
  realizability of a link in a 3-manifold $M$  as a tool for
  understanding the complex structures on the cone $C(M)$.
\end{abstract}
%\keywords{}  
%\bodymatter
\maketitle
\section{Introduction}

%An important corpus of works in singularity theory concentrate on
%answering the following class of questions: which analytical
%properties of a singularity are reflected by its topology?

Let $(Z,p)$ be a normal complex surface singularity, and 
$M_Z$ its link, i.e., the $3$-manifold obtained as the
boundary of a small regular neighbourhood of $p$ in $Z$.  Then, $Z$ is
locally homeomorphic to the cone $C(M_Z)$ on $M_Z$.

Given a surface singularity link $M$, there may exists many different
analytical structures on the cone $C(M)$, i.e., normal surfaces
singularities $(Z,p)$ whose $M_Z$ is homeomorphic to $M$. A natural
problem is to understand these analytic structures on $C(M)$. In this
paper we present an approach by studying the \emph{principal analytic
  link-theory on $M$}. Our aim is to present this point of view to
encourage people to pursue this area.

If $C$ is an analytic curve on $Z$ through $p$, set $L_C = C \cap
M_Z$; the pair $(M_Z,L_C)$, defined up to diffeomorphism, is the {\it
  link} of $C$.  Notice that $L_C$ is a link in $M_Z$ in the usual
topological sense: a disjoint union of circles embedded in a $3$
manifold. In this situation we say $(M,L)=(M_Z,L_C)$ is \emph{
  analytically realized by $(Z,C)$} or simply \emph{analytic}.
%Now, let $M$ be a $3$-dimensional manifold and $L$ a link in $M$. We
%say that $(M,L)$ is {\it analytic} if there exists an analytic
%structure $(Z,p)$ on $C(M)$ and a curve $(C,p)$ on $Z$ such that
%$(M,L)$ is diffeomorphic to the link $(M_Z,L_C)$. 
We say that $(M,L)$
is {\it principal analytic} if it is analytically realized by a pair
$(Z,C)$ such that $C=f^{-1}(0)$, where $f\colon (Z,p) \to (\mathbb
C,0)$ is a germ of holomorphic function. In other words, the ideal of
${\mathcal O}_{Z,p}$ defining $C$ is principal.

% A natural question then arises: given 
Given a $3$-manifold $M$, we then want to %develop the
% analytic (resp.\ principal analytic) link theory on $M$, i.e.,

{\smallskip\par\leftskip\parindent\noindent\hbox to
  0pt{\hss$\bullet$~~}\ignorespaces describe the links $L \subset M$
  which are analytic (resp.\ principal analytic),
\par\leftskip\parindent\noindent\hbox to 0pt{\hss$\bullet$~~}\ignorespaces
describe how these links distribute among the different analytic
  structures on $C(M)$.\par\smallskip}

When $M$ is the $3$ sphere, then $Z$ is smooth \neumanncite{neumann:Mu}, and
one deals with plane curves. The two notions, analytic link and
principal analytic link, then coincide under the classical name
``algebraic link.'' Their classification is one of the main results of
the classical theory of plane curves singularities:  such a
link is built by repeated cabling operations, and a link $L\subset 
  S^3$ is algebraic iff it is obtained by successive cabling
operations which satisfy the so-called Puiseux inequalities.  For
details see \neumanncite{neumann:EN}.
 
For general $M$ the analytic link theory is still well understood, as
a consequence of Grauert \neumanncite{neumann:G}. The main fact is that
analytic realizability is a topological property, so analytic
link-theory is no help to understand the analytic structures on
$C(M)$:
\begin{theorem}\label{neumann:thmanalytic} If $L\subset M$ is analytic,
  it can be realized by a curve $C$ on \emph{any} normal surface
  singularity whose link is $M$.  Moreover, $(M,L)$ is analytic iff it
  can be given by a plumbing graph (see section \ref{neumann:sec2})
  whose intersection matrix is negative definite.
\end{theorem}

In section \ref{neumann:sec2} we show that the existence of
\emph{some} analytic structure on $C(M)$ for which a link is principal
analytic in $M$ is still an easily described topological property (see
Thm.\ \ref{neumann:thmprincipal}, which treats the more general
context of \emph{multilinks} since the zero-set of a holomorphic
function may have multiplicity).  But, in contrast to analytic
realizability, principal analytic realizability depends on the
analytic structure.
% Indeed, if the
% equivalent conditions of Theorem \ref{neumann:thmprincipal} are
% realized and if one fixes the analytic structure $(Z,p)$ on $C(M)$,
% then a holomorphic function $f\colon (Z,p) \to (\mathbb C,0)$ does not
% necessarily exist. 
Here is an explicit example (\neumanncite{neumann:Nemethi}, exercise 2.15):
the two equations $x^2 + y^3 + z^{18} = 0$ and $z^2 + y(x^4 + y^6)=0$
define two germs $(Z,0)$ and $(Z',0)$ in $\mathbb C^3$ with
homeomorphic links (the plumbing graph $\Delta$ is a string o---o---o
of $3$ vertices %$(1), (2)$ and $(3)$
with (genus, Euler number) weights $(1,-1), (0,-2)$ and $(0,-2)$).
But the link of the holomorphic function $f =z: (Z,0) \to (\mathbb
C,0)$ does not have a principal realization in $(Z',0)$.
  Indeed, let $\pi\colon X \to Z$ and $\pi': X'
 \to Z'$ be two resolutions of $Z$ and $Z'$ with dual graph $\Delta$;
 the compact part of the total transform $(f \circ \pi)$ is $E_1 + E_2
 + E_3$ whereas the maximal cycle on $Z'$, realized by the generic
 hyperplane section, is $2E'_1 + 2E'_2 + 2E'_3$, which is greater.

 In section \ref{neumann:sec3} we give some results and conjectures
 about the division of principal analytic links in $M$ among the
 different analytic structures of $C(M)$. In section
 \ref{neumann:sec4} we describe the implications for splice
 singularities, a class of singularities %that are defined 
characterized
 by their principal analytic link theory.

\section{A principal analytic realization theorem}
\label{neumann:sec2}

%First let us recall some facts concerning the plumbing representation
%of the pair $(M_Z,L_C)$.  

We first recall the plumbing representation
of the pair $(M_Z,L_C)$.  

Let $\pi\colon X \to Z$ be a resolution of the normal surface germ
$(Z,p)$ such that $\pi^{-1}(p)$ is a normal crossing divisor with
irreducible components $ E_1,\ldots, E_n$.  Denote by $\Delta$ the
dual graph associated with $\pi^{-1}(p)$, with vertex $v_i$ weighted
$(g_i,e_i)$ by genus $g_i$ and self-intersection $e_i=E_i^2 <0$ of the
corresponding $E_i$ in $X$.
Then $M_Z$ is homeomorphic to the boundary $M$ of the $4$-dimensional
manifold obtained from $\Delta$ by the classical plumbing process
described in \neumanncite{neumann:HNK}.  We call $M$ a \emph{plumbing
  manifold}.

By Grauert \neumanncite{neumann:G}, a plumbing manifold $M$
is the link of a normal complex surface singularity iff it can be
given by a plumbing graph $\Delta$ whose associated intersection
matrix $I(\Delta) = (E_i . E_j)_{1 \leq i,j \leq n}$ is negative
definite.

%Now, let 
If $C$ is a curve on $Z$, and $\pi\colon X \to Z$ a
resolution of $X$ such that $\pi^{-1}(C)$ is a
normal crossing divisor, let $\Delta$ be the dual graph of the divisor
$\pi^{-1}(p)$ decorated with arrows corresponding to the components
%$l_1,\ldots,l_r$ 
of the strict transform of $C$ by $\pi$. Then the
plumbing graph $\Delta$ completely describes the homeomorphism class
of the so-called plumbing link $(M_Z,L_C)$. 

Suppose now that the curve $C$ is the zero locus of a holomorphic
function $f\colon (Z,p) \to (\mathbb C,0)$. Let $C_1,\ldots C_r$ be
the irreducible components of $C$ and let $L_C = K_1 \cup \ldots \cup
K_r$ be its link. Recall that a multilink is a link whose components
are weighted by integers.
%The
%singularity of $f$ at $p$ being allowed to be non isolated, 
One
defines the {\it multilink} of $f$ by $L_f = m_1 K_1 \cup \ldots \cup m_r
K_r$, where $m_j$ is the multiplicity of $f$ along the branch $C_j$.
The plumbing graph $\Delta$ of $(M_Z,L_C)$ is then completed by
weighting each arrow by the corresponding
multiplicity $m_j$.%\comment{Give an example?}

%\begin{definition} 
Let $(M,L)$ be a plumbing multilink with graph
  $\Delta$. For each vertex $v_i$ of $\Delta$, let $b_i$ be the sum of
  the multiplicities carried by the arrows stemming from vertex $v_i$, and
  set $b(\Delta) = (b_1,\ldots,b_n) \in {\mathbb N^n}$. The
  {\it monodromical system} of $\Delta$ is the linear system
$I(\Delta) ^{t}(l_1,\ldots,l_n) + {}^{t}b(\Delta) = 0$  with unknowns
$(l_1,\ldots,l_n)$,  where ${^t}$ means the transposition
%\end{definition}
 \begin{theorem}\label{neumann:thmprincipal} 
Let $L = m_1 K_1 \cup \ldots \cup m_r K_r$ be a
  multilink in a $3$-manifold $M$ with positive multiplicities $m_i$.
  The following are equivalent:
 \begin{description}
%\begin{enumerate}
 \item[(i)] The multilink $(M,L)$ is principal analytically realized
   from some analytic structure $(Z,p)$ on $C(M)$.
 \item[(ii)] $(M,L)$ is a plumbing multilink admitting a plumbing graph
   $\Delta$ whose monodromical system admits a solution
   $(l_1,\ldots,l_n)\in (\mathbb N_{>0})^n$.
 \item[(iii)] $(M,L)$ is analytic and $[m_1 K_1 \cup \ldots
   \cup m_r K_r] = 0$ in $H_1(M,{\mathbb Z})$
 \item[(iv)] $(M,L)$ is a fibered multilink and some power of the
   monodromy $\Phi\colon F \to F$ of the fibration is a product of
   Dehn twists on a collection of disjoint closed curves which
   includes all boundary curves.
\end{description}
% \end{enumerate}
\end{theorem}
\begin{proof} 
  The equivalence between (ii), (iii) and (iv) appears in
  \neumanncite{neumann:EN} when $M$ is a $\mathbb Z$-homology sphere
  and in \neumanncite{neumann:P} (5.4) when $L$ is a link. The methods
  generalize.

  If $f\colon (Z,p) \to (\mathbb C,0)$ is a holomorphic function, then
  a monodromical system is 
nothing but 
%just 
the well known system in
  complex geometry (\neumanncite{neumann:L1}, 2.6): $\big( \forall i =
  1,\ldots,n,\ \ (f).E_i = 0\,\big)$, where $(f)$ is the total
  transform of $f$ in a resolution of $Z$ and $f$ with exceptional
  divisor $E = \sum_{i=1}^n E_i$. Then (i)$\Rightarrow$(ii) is done.
 \par

 (ii)$ \Rightarrow $(i) is proved in \neumanncite{neumann:P} when $L$ is a link
 (5.5). Let generalize the proof to multilinks. The idea is to perform
 a surgery along the multilink $L$ in order to realize the new
 $3$-manifold as the boundary of a degenerating families of curves,
 using a realization theorem of Winters (\neumanncite{neumann:W}).

 Let us consider the plumbing graph $\Delta'$ obtained from $\Delta$
 by replacing 
%the extremity of each arrow $(j)$ 
each arrow $v_i$ o$\!\longrightarrow v_j$
by a string of vertices
 $v_i\hbox{ o---o---o---$\cdots$---o---o}$ as follows: $v_i$ is
 the vertex carrying the arrow $v_j$; set $p_1 = l_i$, $d_1 = m_j$, and
 consider the integers $q_1\geq 1$ and $r_1$ such that $p_1 = q_1 d_1
 - r_1$, with $0\leq r_1 < d_1$. Set $p_2 = d_1$ and $d_2 = r_1$, and
 repeat the process on $p_2$ and $d_2$ by taking $q_2 \geq 1$ and
 $r_2$ such that $p_2 = q_2 d_2 - r_2$, with $0\leq r_2 < d_1$.  Then
 iterate the process until $r_m =0$. The string has $m$
 vertices weighted from $v_i$ by the Euler classes $-q_1, -q_2,
 \ldots, -q_m$ and by genus zero.

 The monodromical system of $\Delta'$ (which has $b(\Delta') = 0$) has
 the following $(l'_k)$ as a solution: $l'_k = l_k$ when $v_k$ is a
 vertex of the subgraph $\Delta$, and $l'_k = d_k$ for the vertex
 $v_k$ of the string carrying $-q_k$. According to
 \neumanncite{neumann:W}, there then exists a degenerating family of
 curves (i.e., a proper holomorphic family which has no critical value
 except $0$) $g:\Sigma \to \{ z \in {\mathbb C} / |z| <1\}$, whose
 special fiber $f^{-1}(0)$ has $\Delta'$ as dual graph
 \neumanncite{neumann:W}.

   By Zariski's lemma (e.g., \neumanncite{neumann:Mu}), condition (ii)
   implies that the intersection matrix $I(\Delta)$ is negative
   definite. One then obtains from $\Sigma$ a normal surface $Z$ by
   shrinking to a single point $p$ all the irreducible components of
   $f^{-1}(0)$ corresponding to vertices of $\Delta$.  Then $g$ induces
   a holomorphic function $Z \setminus \{p\} \to \mathbb C $, which
   extends by $p \mapsto 0$ to a holomorphic function $f\colon Z \to
   \mathbb C$ as $(Z,p)$ is normal. This $f$ realizes $(M,L)$.

\end{proof} 

%\section{The general case: three outlines}
\section{Dependence on analytic structure}
\label{neumann:sec3}
%In this section we point out three partial results --- the
%two first ones through examples --- concerning the distribution of the
%principal analytic realizable links in a $3$-manifold $M$ among the
%different analytical structures on the cone $C(M)$.

%\subsection{On the dependance of the analytic structure}

There exist some $3$-manifolds $M$ whose principal analytic
knot-theory does not depends on the analytic structure.  For example,
when $M$ is the link of a rational singularity, then any principal
analytic realizable multilink $(M,L)$ is so realizable in any analytic
structure $(Z,p)$ on $C(M)$ \neumanncite{neumann:A}. The same conclusion holds
when $M$ is the link of a minimally elliptic singularity and $L$ is a
knot, possibly with multiplicity (\neumanncite{neumann:R}, lemma p.\ 102).  It
seems likely that in most, if not all, other cases the principle
analytic knot theory is sensitive to analytic structure. We are
willing to dare a conjecture in the $\mathbb Z$--homology sphere case.
%\begin{example}  
In this case there is no homological obstruction to the principal
analytic realizability (Theorem \ref{neumann:thmprincipal} (iii)).

Denote the %Pham-
Brieskorn singularity $\{(x_1, x_2,x_3) \in {\mathbb C}^3 ~|~
x_1^p+x_2^q + x_3^{r} = 0\}$ by $V(p,q,r)$ and its link by $M(p,q,r)$.
If $p,q,r$ are pairwise coprime then $M=M(p,q,r)$ is a $\mathbb
Z$--homology sphere.  The only $\mathbb Z$--homology sphere links of
rational and minimally elliptic singularities are $M(2,3,5)$
(rational) and $M(2,3,7)$ and $M(2,3,11)$ (minimally elliptic). So the
principal knot theory of these is completely understood (even the
principal analytic link theory for $M(2,3,5)$).

\begin{conjecture}\label{neumann:conf3.1}
  Let $M$ be a $\mathbb Z$-homology sphere link other than $M(2,3,5),
  M(2,3,7)$, $M(2,3,11)$. Then for any analytic structure $(Z,p)$ on
  $C(M)$ there exists an analytic knot in $M$ which is not realized by
  a holomorphic germ $(Z,p) \to (\mathbb C,0)$.
\end{conjecture}

We can prove the conjecture in many cases. We give two examples to
illustrate the arguments.  $M(p,q,r)$ is Seifert fibered with singular
fibers $L_1,L_2,L_3$ realized as principal analytic knots in
$V(p,q,r)$ by $L_i=M \cap \{x_i = 0\}$.

%\subsection{Incompatible sets of links}
 
1) In $M(2,3,13)$ let $L$ be the $(2,1)$--cable on $L_3$.  As $L$
satisfies condition (iii) of Theorem \ref{neumann:thmprincipal}, it
is principal analytic in some analytic structure.

2) In $M(3,4,19)$ let $L$ be the $(2,3)$--cable on $L_3$.   As $L$
satisfies condition (iii) of Theorem \ref{neumann:thmprincipal}, it
is principal analytic in some analytic structure.
\begin{proposition} \label{neumann:propex}
  (1) Let $(Z,p)$ be an analytic structure on $C(M(2,3,13))$ such that
  $L_3$ is realized by a holomorphic function $f_3: (Z,p) \to
  (\mathbb C,0)$. Then $L$ is not realized by any $f\colon (Z,p) \to
  (\mathbb C,0)$ on $(Z,p)$.

  (2) Let $(Z,p)$ be an analytic structure on $C(M(3,4,19))$ such that
  both $L_2$ and $L_3$ are realized on $(Z,p)$ by $f_2,f_3\colon (Z,p) \to
  (\mathbb C,0)$  respectively.
  Then, $L$ is not realized by any $f\colon (Z,p) \to (\mathbb
  C,0)$.
\end{proposition}
\begin{proof}
  (1) Assume the contrary.  Let $E_i, i=1,\ldots,5$ be the irreducible
  components of the exceptional divisor of the minimal resolution
  $\pi\colon \Sigma \to Z$ of $(Z,p)$ as in the figure below.

 The total transform of $f_3$ is $(f_3 \circ \pi) = 3 E_1 + 2 E_2 + 6
 E_3 + E_4 + E_5 + l_3$, so its multiplicity on
  $E_4$ is $1$, which means that $f_3$ is a local coordinate on
  the transverse curve $f^{-1}(0)$ to $E_4$, so $f^{-1}(0)$ is smooth.
  Therefore, the Milnor fibre $F_t$ of $f$, $t \neq 0$, is a disk, as
  it is a smoothing of $f^{-1}(0)$. But the total transform of $f$ is
  $(f \circ \pi) = 6 E_1 + 4 E_2 + 12 E_3 + 2E_4 + E_5 + l$, which
  leads to $\chi(F_t) = 6 + 4 -12 -2 + 1= -3$. Contradiction.
$$
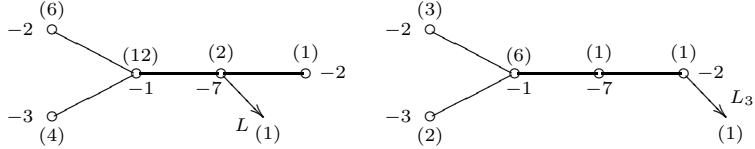
\xymatrix@R=12pt@C=15pt@M=0pt@W=0pt@H=0pt{
\leftNeTag{\overNeTag\CirC{(6)}{9pt}}{-2}{6pt}\ar@{-}[drr]
&&&&&&&&&
\leftNeTag{\overNeTag\CirC{(3)}{9pt}}{-2}{6pt}\ar@{-}[drr]\\
&&\CirC\ar@{-}[rr]^(.05){(12)}_(.05){-1}_(.85){-7}&&
\overNeTag\CirC{(2)}{9pt}\ar@{-}[rr]\ar[dr]_(.85){L}&&
\rightNeTag{\overNeTag\CirC{(1)}{9pt}}{-2}{6pt}&&&
&&\CirC\ar@{-}[rr]^(.05){(6)}_(.05){-1}_(1){-7}&&
\overNeTag\CirC{(1)}{9pt}\ar@{-}[rr]&&
\rightNeTag{\overNeTag\CirC{(1)}{9pt}}{-2}{6pt}\ar[dr]^(.9){L_3}
\\
\leftNeTag{\underNeTag\CirC{(4)}{2pt}}{-3}{6pt}\ar@{-}[urr]
&&&&&\underNeTag{}{~(1)}{3pt}&&&&
\leftNeTag{\underNeTag\CirC{(2)}{2pt}}{-3}{6pt}\ar@{-}[urr]
&&&&&&&\underNeTag{}{~(1)}{3pt}\\&
}
$$

2) Assume the contrary to the proposition. Let $f\colon (Z,p) \to
(\mathbb C,0)$ be such that $L$ is the link of $f^{-1}(0)$. The splice
diagram is as follows:
$$
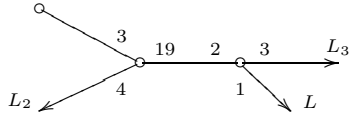
\xymatrix@R=16pt@C=18pt@M=0pt@W=0pt@H=0pt{
\CirC\ar@{-}[drr]^(.75)3\\
&&\CirC\ar@{-}[rr]^(.25){19}^(.75)2\ar[dll]^(.25)4_(1.1){L_2}&&
\CirC\ar[rr]^(.25)3^(1){L_3}\ar[dr]_(.25)1^(1.1)L&&\\
&&&&&\\
}
$$
The semi-group of values $\Gamma(C)$ of the curve $C = f^{-1}(0)$
contains the two values $2$ and $9$ associated with the functions
$f_2$ and $f_3$ (these are computed as the product of splice diagram
weights adjacent to the path between the corresponding arrowheads).
%(\neumanncite{neumann:C}, section 2). 
So $\Gamma(C)$ contains the
semi-group
$\langle 2,9\rangle = \{ 2,4,6,8, 8 + {\mathbb Z}^{+} \ldots \}$
which has $4$ missing numbers ($1,3,5$ and $7$). Therefore, the
$\delta$-invariant $\delta(C)$ of $C$, which counts the number of gaps
in $\Gamma(C)$, %(see \neumanncite{neumann:C}, section 2) 
is at most $4$, so
$\mu(C) = 2 \delta \leq 8$.
But the multiplicities of $C$ leads to $\chi(F_t) = 12+9-36 -6 + 2 =
-19$. Then $\mu = 20$. Contradiction.
\end{proof} 
 
%\subsection{Ubiquitous links}

In view of the above discussion, the following question is natural.
\begin{question}
Let $M$ be a surface singularity link. Do there exist
analytic links in $M$ that have the {\it ubiquity property}, i.e.,
that are principal analytic in any analytic structure $(Z,p)$ with
link $M$?  
\end{question}

\iffalse
A positive answer is given by:

\begin{theorem}
  {\bf (\neumanncite{neumann:CNP}, 4.1.)} Let $(Z,0)$ be any analytic
  structure on $C(M)$, and let $\pi\colon (\Sigma,E) \to (Z,0)$ be a
  good resolution with $E = \sum_{i=1}^n E_i$. For any purely
  exceptional divisor $D = \sum m_i E_i$ such that
 $$\forall i \in \{1,\ldots,n\}, (D + E + K_{Z}). E_i \leq 0,$$
 there exists a holomorphic germ $f\colon (Z,p) \to (C,0)$ whose total
 transform by $\pi$ is the divisor $(f \circ \pi) = D + Strict (f \circ
 \pi) $, its strict part consisting of $n_i \geq 2 g(E_i) + 1$ smooth
 components transversal to $E_i$ at smooth point of $E$ and one smooth
 transversal component through each intersection point $E_i \cap E_j,
 i \neq j$.
\end{theorem}
As the configuration of $E$ and the numerical data $D = \sum m_i E_i$
do not depend on the analytic structure $(Z,p)$, the the link $L_f$
has the ubiquity property. This link has very many components. We can
reduce the number somewhat, but ``ubiquitous links'' probably always
need many components.
\fi
 A positive answer is given by the following theorem, which is a
 slight improvement of Theorem 4.1 of \neumanncite{neumann:CNP}:

\begin{theorem}
  Let $g_i$ be the genus weights of a negative definite plumbing for
  $M$ and $L$ be the $\left(\sum (2g_i+1)\right)$--component link
  consisting of the boundaries of $2g_i+1$ transverse disks to each
  base curve $E_i$ of the plumbing. Then $L$ is principal analytic in
  any analytic structure $(Z,0)$ on $C(M)$.
\end{theorem}
If $M$ is a $\mathbb Q$--homology sphere this ubiquitous link has one
component for each vertex of the negative definite plumbing graph
(since each $g_i$ is $0$). It seems unlikely that this can be reduced
much in general, although in the rational and minimally elliptic cases
$1$ component suffices.

In the original version of this theorem in \neumanncite{neumann:CNP}, L
had one extra component for each edge of the plumbing graph. The
proof of the above improvement is a small modification of that proof,
but the 8 page limit on contributions to this volume leaves
insufficient space to give it here.

\section{Singularities of splice type}
 \label{neumann:sec4}

 In \neumanncite{neumann:splice, neumann:zsphere} the first author and
 J. Wahl introduce an important class of singularities with $\mathbb
 Q$--homology sphere links called ``splice-quotient singularities,''
 or briefly ``splice singularities.'' This class includes all rational
 singularities and all minimal elliptic singularities with $\mathbb
 Q$--homology sphere links. This was first proved by Okuma
 \neumanncite{neumann:O}, see also the appendix to
 \neumanncite{neumann:splice}. It now follows from the following
 characterization of splice singularities by their principal analytic
 knot theory. We assume $M$ is a $\mathbb Q$--homology sphere. The
 dual resolution graph is then a tree.

\begin{theorem}
[End Curve Theorem (\neumanncite{neumann:endcurves}]
To each leaf of the dual resolution graph there
is associated a knot in $M$, namely the link of a transverse curve to
the corresponding exceptional curve of the resolution; the
  singularity is splice if and only if each of these knots, taken with
  some multiplicity, has a principal analytic realization.
\end{theorem}

The property of being rational or minimally elliptic is topologically
determined, so any singularity with one of these topologies is splice.
But in general the same topology may support both splice and
non-splice singularities. For example, 
%the Brieskorn singularities are splice, 
Proposition \ref{neumann:propex} implies for either of the examples it
addresses, that if one takes an analytic
structure for which the knot $L$ is principal, then that singularity
is not splice.

Associated to any plumbed $\mathbb Q$--homology sphere is a simplified
version of the resolution graph, called the splice diagram. A
fundamental property of the splice diagrams of links of splice
singularities is the so-called ``semigroup condition'' (loc.\ cit.).
The following question is of fundamental importance, since a $\mathbb
Z$--sphere counterexample would give a complete intersection
singularity with $\mathbb Z$-homology sphere link that is not splice
(conjectured not to exist), and would give a likely candidate also to
contradict the Casson Invariant Conjecture of
\neumanncite{neumann:casson}.

\begin{question}
  Is the principal knot theory of a splice singularity 
%with $\mathbb  Z$-homology sphere link 
topologically
  determined? Specifically, is a knot $L\subset M$ with multiplicity
  principal analytic if and only if it represents zero in homology of
  $M$ and the splice diagram for $(M,L)$ satisfies the semigroup
  condition?
\end{question}

\smallskip\noindent{\bf Acknowledgments:} The first author
acknowledges the support of the NSF and the NSA for this research.

\end{document}